\documentclass[a4paper,11pt]{article}
\usepackage{palatino}
\usepackage{booktabs}
\usepackage{url}
\usepackage[whole]{bxcjkjatype}
\usepackage{vmargin}
\setmarginsrb{2cm}{2cm}{1.95cm}{1.95cm}{0mm}{0mm}{0mm}{8mm}
\usepackage{enumitem}
\usepackage{multirow}
\usepackage{authblk}
\usepackage{marvosym}
\usepackage{graphicx}
\usepackage{cite}
\usepackage[colorlinks=true, allcolors=blue, breaklinks=true]{hyperref}

\title{\vspace*{-1.5cm} \bfseries Revisiting Javanese \emph{pranata mangsa}:\\ On ethnic groups and the four sample cities in Java}
\author{\normalsize Natanael Karjanto\thanks{\Letter: \url{natanael@skku.edu} \href{https://orcid.org/0000-0002-6859-447X}{\includegraphics[scale=0.08]{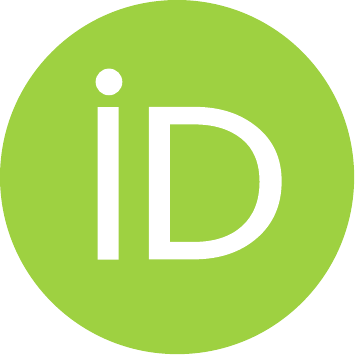}}}}
\affil{Department of Mathematics, University College, Natural Science Campus\\ Sungkyunkwan University, Suwon~16419, Republic of Korea}

\date{\vspace*{-0.5cm} \scriptsize Updated \today}

\begin{document}
\maketitle \vspace*{-1cm}
\begin{abstract}
\emph{Pranata mangsa} is a local knowledge based on a calendrical system utilized historically by Javanese peasant farmers in conducting their agricultural activities in cultivating plants and vegetation. It is also used by fishermen to guide them not only in capturing fish and other aquatic animals but also in predicting the type of seafood they might gather. Although the Javanese possess and maintain their own calendrical system based on a combination of solar and lunar calendars, \emph{pranata mangsa} is solely based on the solar calendar and passed on from one generation to the next verbally. Literature has confirmed that by combining this local knowledge and scientific data, the Javanese community has a better resilience in adapting to extreme hydrological events due to global warming and climate change. In this comment, we provide remarks on ethnic groups and the four sample cities in Java where the latter were selected for the \emph{pranata mangsa} study.
\end{abstract}

\section{Introduction}

Global warming and climate change play an influential impact on the agricultural sector, particularly in developing countries~\cite{nelson2009climate}. It is essential for the government, policy makers, and farmers to implement mitigation strategies and acquire sustainable solutions in the onset of extreme events~\cite{arora2019impact,malhi2021impact}. One of such scenarios is by employing local knowledge. Among the Javanese farmers in Indonesia, a solar-based agricultural calendar called \emph{pranata mangsa} has been used for over a millennia~\cite{daldjoeni1984pranatamangsa}. In particular, Zaki et al. provide excellent coverage on the local knowledge of \emph{pranata mangsa} in the Javanese society with the attempts among local farmers in adjusting to various extreme hydrological events, including flood, drought, and harvest failure~\cite{zaki2020adaptation}. This comment covers several remarks on the paper.

\section{Comment}
\subsection{Ethnic groups in Java}
According to Statistics Indonesia, the population of Java in 2017 is 148,173,100. A recent census in 2020 shows that there are 151,900,00 people residing in Java, around 56.10\% of the total Indonesian population. Zaki et al.'s estimate of 150~million people at the time of their writing is an excellent one~\cite{zaki2020adaptation}. However, the authors mentioned that ``there are three ethnic groups in Java: Betawi, Sundanese, and Javanese.'' This is incomplete at best and ignorance at worst. Even in the absence of the word ``only'', their sentence is still inaccurate. Adding the phrase ``at least'' might soften the error, but the list that follows is still deficient. We would like to clarify that the predominant ethnic groups in Java according to their percentage are Javanese, Sundanese, Madurese, Betawi, and Chinese. 

Certainly, there are many other ethnic groups from other parts of Indonesia residing in Java, as well as some as a result of a mixed marriage between Indonesians and non-Indonesians. Indeed, only Javanese and Sundanese are native to the island. The third group Madurese has been immigrating to the eastern part of Java since the 18th century from their native island of Madura, off the northeast coast of Java~\cite{hefner1997java}. The fourth group Betawis are an Austronesian ethnic group ``native'' to Batavia, the colonial term for Jakarta, and its immediate surroundings. Furthermore, since the Betawi people emerged in the 18th century as an amalgamation of various immigrant ethnic groups into Batavia, the adjective ``native'' with the quotation mark makes more sense. 

The fifth group had been arriving from China to Indonesia in several stages since the 13th century. The majority came from the southern part of China, which now constitutes the Fujian and Guangdong provinces in modern China. Some Chinese Indonesians hesitate to reveal their ethnic identity due to fears of stigmatization and discrimination, and thus it might be hard to obtain their demographic information~\cite{suryadinata1976indonesian}. However, a comprehensive study from 2016 revealed that there are around 1.2\% of Chinese Indonesians residing in the country, contrary to other sources that mention in the 3\%--3.3\% range~\cite{mackie2005how,arifin2017chinese,chinesedisaspora,facts2015chinese,databoks2016indonesia,minority2018indonesia}. See also~\cite{lindsey2005chinese,turner2007chinese,suryadinata2008chinese,dieleman2011chinese,setiadji2017chinese,rakhmat2021indonesia} for a well-rounded perspective on Chinese Indonesians. 

In the authors' list, leaving out Chinese Indonesians can be understandable, as historically, they did not belong to the indigenous ethnic group in Indonesia. But we should not exclude the Madurese from the list. If the authors meant ethnic groups native to the island of Java, then Betawi people should be left out as well. When the \emph{pranata mangsa} was established circa the 11th century, it was most likely the Chinese had not arrived in Java yet, and although we could not prove the absence of the Madurese people during that time, the Betawis were certainly nonexistent.

Many Chinese Indonesians, particularly from Java, have made notable contributions in various areas, including but not limited to art, academics, business, politics, and sports. Some examples in what follows are by no means exhaustive. In music, we have traditional Chinese musician Inke Kirei (Dian Prasidha) and violinist Kezia Amelia Angkadjaja from Surabaya. In sports, particularly badminton, we have the couple Susi Susanti (王蓮香) from Tasikmalaya and Alan Budikusuma (魏仁芳) from Surabaya, both won gold medals in the 1992 Olympic Games, as well as Rudy Hartono (梁海量), also from Surabaya, an 8-time winner of the All-England Cup, an unbreakable record up to this date. 

Some academicians tend to be less well-known, but they have significant contributions to their respective areas of expertise. We have a mathematics educator Koko Martono from Garut (1952--2018), mathematicians Hendra Gunawan from Bandung, and Marcus Wono Setya Budhi from Weleri, Kendal. There are also physicists Tjia May On from Banyuwangi (1934--2019), The Houw Liong from Bandung, and Tak Ik Gie (\quad--2004). Chinese Indonesians diaspora overseas include multi-disciplinarian Ken Kawan Soetanto (陳文權) from Surabaya (Waseda University, Japan), immunochemist Yow-Pin Lim from Cirebon (Brown University, Providence, Rhode Island), scholar Merlyna Lim from Dayeuhkolot (Carleton University, Canada), applied mathematician Djoko Wirosoetisno from Surabaya (Durham University, England), number theorist Ade Irma ``Chacha'' Suriajaya (謝佳汶) from Jakarta (Kyushu University, Japan), Hartono Tjoe (Pennsylvania State University), and Melisa Hendrata (California State University Los Angeles).

\subsection{The four cities climate comparison}

When discussing the study area of local knowledge in the Javanese society, the authors have selected particular four cities in Java in their sample studies. One town in West Java (Indramayu), two in Central Java (Sleman and Sukoharjo), and one in East Java (Ngawi). Since there was no given explanation, many readers wonder why the authors have chosen those particular four cities. The best clue might be gleaned from Figure~3 of the paper, where it displays histograms of the precipitation rate for the entire (annual) \emph{pranata mangsa} period, where the data from the National Aeronautics and Space Administration (NASA) and Japan Aerospace Exploration Agency (JAXA) were gathered for the period of 18~years (1998--20215). From these histograms, it seems that the four cities have similar precipitation rates and only exhibit a small variation annually. 

According to the Köppen--Geiger climate classification system, although predominantly belonging to Group A (tropical climates), the island of Java and its larger provinces themselves have several climate variations. In West Java, the most prevalent ones are \emph{Af} and \emph{Am} (tropical rainforest and tropical monsoon climates, respectively), whereas Central and East Java are dominated by \emph{Am} and \emph{Aw}, respectively. \emph{Aw} means tropical savanna climate with dry-winter characteristics. While the average precipitation for \emph{Af} is at least 60~mm monthly, in \emph{Am}, \emph{Aw}, or \emph{As}, a similar rate occurs during the driest months only~\cite{peel2007updated,beck2018present,beck2020present}.

Looking at these four cities in more detail, we discover some interesting facts. Indramayu in West Java has a tropical savanna climate (\emph{Aw}) with moderate to little rainfall from May to November and heavy rainfall from December to April. The average rainfall is around 1700~mm annually. Sleman is actually located in the Special Administrative Region of Yogyakarta instead of Central Java, albeit its northern border shares with the latter. It features a tropical monsoon climate (\emph{Am}), with a long wet season from October until June the following year and a shorter dry season from July until September. Its annual precipitation rate is around 2200~mm. Sukoharjo is a regency in Central Java that also has a tropical monsoon climate (\emph{Am}), with little to moderate rainfall from June to October and a heavier one from November until May. It also boasts an annual precipitation rate above 2000~mm. Ngawi is located in the western part of East Java, sharing a border with Central Java, midway from north to south. Similar to the previous two cities, Ngawi is also characterized by a tropical monsoon climate (\emph{Am}), with heavy rains at the beginning and end of the dry season in the middle of the year. With its precipitation rate of 1600~mm annually, it is the driest city among the four. See Table~\ref{table} for a climate comparison of these four sample cities.
\begin{table}[h] \vspace*{-0.5cm}
\caption{A comparison of the four sample cities. We observe that Indramayu has a different climate from the other three cities due to its geographical location as a coastal city, whereas the other three sample cities were located relatively far from the sea.}	\label{table}
\begin{center}
\begin{tabular}{@{}llccc@{}}
\toprule	
	 		&  					&  				&  	 		& Annual \\ 	
City		& Province			& Elevation~(m)	& Climate	& precipitation \\ 
			&					&				&			 &  rate~(mm)	\\				\midrule
Indramayu 	& West Java 		& 3				& \emph{Aw}  & 1700 \\ 
\multirow{4}{*}{Sleman}		& Yogyakarta  		& \multirow{4}{*}{100--2500} 	& \multirow{4}{*}{\emph{Am}}  & \multirow{4}{*}{2200} \\
            & Special			& 				&			 & 		\\	
			& Administrative 	&				&			 &		\\
			& Region			&				&			 &		\\
Sukoharjo   & Central Java 		& 80--125		& \emph{Am}  & 2100 \\
Ngawi		& East Java 		& 47--500		& \emph{Am}	 & 1600 \\
\bottomrule	
\end{tabular}
\end{center}
\end{table}

Historically, the classical version of the \emph{pranata mangsa} applies for the region between Mount Merapi, an active stratovolcano located on the border between the Special Region of Yogyakarta and the province of Central Java ($7^\circ 32' 26.99''$S, $110^\circ 26' 41.34''$E), and Mount Lawu, also a composite volcano straddling the border between Central and East Java provinces ($7^\circ 37' 30''$S, $111^\circ 11' 30''$E)~\cite{daldjoeni1984pranatamangsa}. Indeed, both Sleman and Sukoharjo are located between these two mountains, and Ngawi is around 50~km northeast direction from Mount Lawu, whereas Indramayu is nearly 400~km northwest direction from Mount Merapi. It would be interesting to investigate whether other areas between these two mountains, as well as other places in Java with similar climatic features to this southeast region of Central Java, exhibit robust characteristics of \emph{pranata mangsa}. The former includes Boyolali, Klaten, Surakarta (Solo), Karanganyar, and the northern portion of Wonogiri. The latter might comprise Cirebon, Majalengka, Kuningan, or North Sumedang in West Java, as well as Blitar, Malang, Madiun, Trenggalek, or Bondowoso in East Java.  

\subsection{Cicadidae}
Zaki et al. mentioned that the insect family Cicadidae, the true cicadas, features the rainy season \emph{rendheng} with its appealing sound. The cicada species that the Javanese farmers encounter might be \emph{Tacua speciosa} or \emph{Purana tigrina} group, which are commonly found not only in Java but also in other parts of Southeast Asia~\cite{duffels2007revision}. In Indonesia, cicadas' sound usually appears at the end of the wet season, a clear sign that the dry season is around the corner and the rain would stop. Table~2 in their paper lists cicadas' acoustic signals in March, whereas the transitional season from the wet to dry seasons is usually around the end of March and April, cf.~\cite{sudibyakto2018manajemen}. We would not really be sure whether the shifting appearance of cicadas is due to global warming and climate change. The global temperature on Earth has been steadily increasing since the industrial revolution at the end of the 19th century. 

In subtropical regions like in Korea, cicadas mostly emerge in late June through August annually. And although there are nearly 3400 species of cicadas exist worldwide, periodical cicadas that emerge altogether once every 13 or 17 years are only found in the eastern part of the United States. They usually appear from late April through early May~\cite{wong2021broox,kritsky2021one}. The connection between the time of the appearance of cicadas in different parts of the world and whether climate change contributes to the shifting in their annual emergence might be a topic of interest for further investigation among colleagues who specialize in insect behavior, ecology, and climate change. See also~\cite{taylor2021cicadian} for the impact on human health.

\section{Conclusion}

We have briefly discussed topics related to the Javanese local knowledge of \emph{pranata mangsa} in this comment. In particular, we have rectified the list of ethnic groups in Java, which should include Madurese and the often neglected, forgotten, and stigmatized Chinese Indonesians. We have also provided a perspective on the four sample cities considered in the study of the \emph{pranata mangsa}, with a proposal of considering other regions with similar climate characteristics to the region where the classical \emph{pranata mangsa} was initially inaugurated. The appearance of cicadas seems to vary not only depending on the geographical location but also due to the potential influence of global warming and climate change. We hope that this comment will stimulate further discussion and future collaboration.


{\small
\begin{thebibliography}{99}
\bibitem{nelson2009climate} Nelson, G. C., Rosegrant, M. W., Koo, J., Robertson, R., Sulser, T., Zhu, T., ... and Lee, D. (2009). \textit{Climate Change: Impact on Agriculture and Costs of Adaptation}. Washington, DC, US: International Food Policy Research Institute.

\bibitem{arora2019impact} Arora, N. K. (2019). Impact of climate change on agriculture production and its sustainable solutions. \textit{Environmental Sustainability}, 2(2), 95--96.

\bibitem{malhi2021impact} Malhi, G. S., Kaur, M., and Kaushik, P. (2021). Impact of climate change on agriculture and its mitigation strategies: A review. \textit{Sustainability}, 13(3), 1318.
	
\bibitem{daldjoeni1984pranatamangsa} Daldjoeni, N. (1984) Pranatamangsa, the Javanese agricultural calendar--Its bioclimatological and sociocultural function in developing rural life. \textit{Environmentalist}, 4(7),15--18.
	
\bibitem{zaki2020adaptation} Zaki, M. K., Noda, K., Ito, K., Komariah, K., Sumani, S., and Senge, M. (2020). Adaptation to extreme hydrological events by javanese society through local knowledge. \textit{Sustainability}, 12(24), 10373.	

\bibitem{bps1} Bekasi Municipality Statistics. Retrieved from  \url{https://bekasikab.bps.go.id/statictable/2021/07/24/2994/jumlah-penduduk-menurut-provinsi-di-indonesia-ribu-2013-2017.html}. Last accessed \today.

\bibitem{bps2} Statistics Indonesia. Data released on 21 January 2021. Retrieved from  \url{https://www.bps.go.id/publication/2021/01/21/213995c881428fef20a18226/potret-sensus-penduduk-2020-menuju-satu-data-kependudukan-indonesia.html}. Last accessed \today.

\bibitem{hefner1997java} Hefner, R. (1997). \textit{Java}. Singapore: Periplus Editions.

\bibitem{knorr2014creole} Knorr, J. (2014). \textit{Creole Identity in Postcolonial Indonesia}. Volume 9 of Integration and Conflict Studies. Berghahn Books.

\bibitem{suryadinata1976indonesian} Suryadinata, L. (1976). Indonesian policies toward the Chinese minority under the new order. Asian Survey, 16(8), 770-787.

\bibitem{mackie2005how} Mackie, J. (2005). How many Chinese Indonesians?. \textit{Bulletin of Indonesian Economic Studies}, 41(1), 97--101.

\bibitem{arifin2017chinese} Arifin, E. N., Hasbullah, M. S., and Pramono, A. (2017). Chinese Indonesians: how many, who and where?. \textit{Asian Ethnicity}, 18(3), 310-329.

\bibitem{chinesedisaspora} n.d. Chinese diaspora across the world: a general overview. \textit{Academy for Cultural Diplomacy}. Retrieved from  \url{https://www.culturaldiplomacy.org/academy/index.php?chinese-diaspora}. Last accessed \today.

\bibitem{facts2015chinese} n.a. (2015). Chinese in Indonesia. \textit{Facts and Details}. Retrieved from \url{https://factsanddetails.com/indonesia/Minorities_and_Regions/sub6_3a/entry-3993.html}. Last accessed \today.

\bibitem{databoks2016indonesia} n.a. (2016). Indonesia, populasi etnis Cina terbanyak di dunia. (Indonesia, the largest ethnic Chinese population in the world.) \textit{Databoks}. Retrieved from  \url{https://databoks.katadata.co.id/datapublish/2016/12/13/indonesia-populasi-etnis-cina-terbanyak-di-dunia}. Last accessed \today, (in Indonesian).

\bibitem{minority2018indonesia} n.a. (2018). Indonesia--Chinese. \textit{Minority Rights}. Retrieved from \url{https://minorityrights.org/minorities/chinese-3/}. Last accessed \today.

\bibitem{lindsey2005chinese} Lindsey, T., and Pausacker, H. (Eds.). (2005). \textit{Chinese Indonesians: Remembering, Distorting, Forgetting}. Singapore: Institute of Southeast Asian Studies.

\bibitem{turner2007chinese} Turner, S., and Allen, P. (2007). Chinese Indonesians in a rapidly changing nation: Pressures of ethnicity and identity. \textit{Asia Pacific Viewpoint}, 48(1), 112--127.

\bibitem{suryadinata2008chinese} Suryadinata, L. (Ed.) (2008). \textit{Ethnic Chinese in Contemporary Indonesia}. Singapore: Institute of Southeast Asian Studies.

\bibitem{dieleman2011chinese} Dieleman, M., Koning, J., and Post, P. (Eds.). (2011). \textit{Chinese Indonesians and Regime Change}. Leiden, the Netherlands: Brill.

\bibitem{setiadji2017chinese} Setijadi, C. (2017). Chinese Indonesians in the eyes of the \emph{pribumi} public. \textit{Institute of Southeast Asian Studies Yusof Ishak Institute Perspective}, 2017(73), 1--12. Retrieved from \url{https://think-asia.org/bitstream/handle/11540/7545/ISEAS_Perspective_2017_73.pdf?sequence=1}. Last accessed \today.

\bibitem{rakhmat2021indonesia} Rakhmat, M. Z. (2021). Indonesia tries to embrace Chinese language but problems persist. \textit{The Conversation}. Retrieved from \url{https://theconversation.com/indonesia-tries-to-embrace-chinese-language-but-problems-persist-165608}. Last accessed \today.

\bibitem{britannica} n.a. (n.d.) Chinese and other Indonesian peoples. \textit{Encyclopedia Britannica}. Retrieved from \url{https://www.britannica.com/place/Indonesia/Chinese-and-other-Indonesian-peoples}. Last accessed \today.

\bibitem{peel2007updated} Peel, M. C., Finlayson B. L., and McMahon, T. A. (2007). Updated world map of the Köppen--Geiger climate classification. \textit{Hydrology and Earth System Sciences}, 11(5): 1633--1644. 

\bibitem{beck2018present} Beck, H. E., Zimmermann, N. E., McVicar, T. R., Vergopolan, N., Berg, A., Wood, E. F. (2018). Present and future Köppen-Geiger climate classification maps at 1-km resolution. \textit{Scientific Data}, 5: 180214.

\bibitem{beck2020present} Beck, H. E., Zimmermann, N. E., McVicar, T. R., Vergopolan, N., Berg, A., Wood, E. F. (2020). Publisher correction: Present and future Köppen-Geiger climate classification maps at 1-km resolution. \textit{Scientific Data}, 7: 274.

\bibitem{duffels2007revision} Duffels, J. P., Schouten, M. A., and Lammertink, M. (2007). A revision of the cicadas of the Purana tigrina group (Hemiptera, Cicadidae) in Sundaland. \textit{Tijdschrift voor Entomologie}, 150(2), 367.

\bibitem{sudibyakto2018manajemen} Sudibyakto, H. A. (2018). Manajemen bencana di Indonesia ke mana? Yogyakarta, Indonesia: UGM Press.

\bibitem{wong2021broox} Wong, K. and Sinnen, C. (2021). Brood X cicadas are emerging at last. \textit{Scientific American}, 324(6), 54--59. 
	
\bibitem{kritsky2021one} Kritsky, G. (2021). One for the books: The 2021 emergence of the periodical cicada brood X. \textit{American Entomologist}, 67(4), 40--46.

\bibitem{taylor2021cicadian} Taylor, C. A. (2021). Cicadian rhythm: Insecticides, infant health and long-term outcomes. \textit{Center for Environmental Economics and Policy} Working Paper No.~9, 45~pp. Retrieved from \url{https://ceep.columbia.edu/sites/default/files/content/papers/n9_v2.pdf}. Last accessed \today.

\end{thebibliography}
}
\end{document}